\documentclass[12pt,a4paper]{article}
\usepackage[cp1251]{inputenc}
\usepackage[T2A]{fontenc}
\usepackage[english]{babel}

\usepackage{amsmath,amsfonts,amssymb,mathrsfs,amsopn,indentfirst,amscd}

\usepackage{graphicx}
\usepackage{amsthm}
\usepackage{hyperref}

\usepackage{amssymb}
\usepackage{amsmath}

\newtheorem*{corollary*}{Corollary}

\begin{document}

\title{The Minimum Norm of a Projector \\ under Linear Interpolation on a Euclidean Ball}

\author{Mikhail Nevskii\footnote{ Department of Mathematics,  P.G.~Demidov Yaroslavl State University,\newline Sovetskaya str., 14, Yaroslavl, 150003, Russia, \newline
               mnevsk55@yandex.ru,  orcid.org/0000-0002-6392-7618 } 
               }
      
\date{May 1, 2023}
\maketitle

\begin{abstract}

\smallskip
We prove the following proposition. 
Under linear interpolation on a Euclidean $n$-dimensional ball $B$, an interpolation 
projector whose nodes coincide with the vertices of a regular simplex inscribed into the boundary sphere has the minimum $C$-norm. This minimum norm $\theta_n(B)$ $=$
 $\max\{\psi(a_n),\psi(a_n+~1)\}$, where
$\psi(t)=\dfrac{2\sqrt{n}}{n+1}\Bigl(t(n+1-t)\Bigr)^{1/2}+
\left|1-\dfrac{2t}{n+1}\right|$, \ $0\leq t\leq n+1$,  \ and~$a_n=\left\lfloor\dfrac{n+1}{2}-\dfrac{\sqrt{n+1}}{2}\right\rfloor$. For any $n$,
$\sqrt{n}\leq \theta_n(B)\leq \sqrt{n+1}.$
Moreover, $\theta_n(B)$ $=$ $\sqrt{n}$ only for~$n=1$ and
$\theta_n(B)=\sqrt{n+1}$ \  iff \ 
$\sqrt{n+1}$  is an integer.

\smallskip 



\end{abstract}

\section{Main definitions}\label{nev_s1}
Further $n\in{\mathbb N}$. An element $x\in {\mathbb R}^n$ is written in the form $ x=(x_1,\ldots,x_n)$. Let $K$ be a convex body in
${\mathbb R}^n$, i.\,e., a compact convex subset of ${\mathbb R}^n$ with nonempty interior. Denote by   
$C(K)$ a space of continuous functions 
$f:K\to{\mathbb R}$ with  the~uniform norm
$$\|f\|_{C(K)}:=\max\limits_{x\in K}|f(x)|.$$
By $\Pi_1\left({\mathbb R}^n\right)$ we mean a set of polynomials in
$n$ variables of degree
$\leq 1$, i.\,e.,  a~set \linebreak of  linear functions on ${\mathbb R}^n$.
For $x^{(0)}\in {\mathbb R}^n,  R>0$, by $B(x^{(0)};R)$ we denote \linebreak  the $n$-dimensional Euclidean ball given
by the inequality 
$\|x-x^{(0)}\|\leq R$. Here 
$$\|x\|:=\sqrt{\langle x,x\rangle}=\left(\sum\limits_{i=1}^n x_i^2\right)^{1/2}.$$
By definition,  $B_n:=B(0;1)$. 

Let $S$ be a nondegenerate simplex in  ${\mathbb R}^n$ with the vertices 
$x^{(j)}=\left(x_1^{(j)},\ldots,x_n^{(j)}\right),$ 
$1\leq j\leq n+1.$ Consider {\it the vertex matrix} of this simplex:
$${\bf S} :=
\left( \begin{array}{cccc}
x_1^{(1)}&\ldots&x_n^{(1)}&1\\
x_1^{(2)}&\ldots&x_n^{(2)}&1\\
\vdots&\vdots&\vdots&\vdots\\
x_1^{(n+1)}&\ldots&x_n^{(n+1)}&1\\
\end{array}
\right).$$
Assume ${\bf S}^{-1}$ $=(l_{ij})$. Linear  polynomials 
$\lambda_j(x):=
l_{1j}x_1+\ldots+
l_{nj}x_n+l_{n+1,j}$
have the~property  
$\lambda_j\left(x^{(k)}\right)$ $=$ 
$\delta_j^k$.
We call $\lambda_j$ {\it the basic Lagrange polynomials corresponding to $S$.}
For~an~arbitrary $x\in{\mathbb R}^n$, 
$$x=\sum_{j=1}^{n+1} \lambda_j(x)x^{(j)}, \quad \sum_{j=1}^{n+1} \lambda_j(x)=1.$$
This means that $\lambda_j(x)$ are {\it  the barycentric coordinates of $x$}
(see, e.\,g.,~\cite{nevskii_1}--~\cite{nevskii_11}).

We say that the interpolation projector
 $P:C(K)\to \Pi_1({\mathbb R}^n)$ corresponds to~a~simplex $S\subset K$  if the nodes of $P$
 coincide with the vertices of $S$. This projector is defined by~the~equalities
$Pf\left(x^{(j)}\right)=
f\left(x^{(j)}\right).$
The following analogue of the Lagrange interpolation formula holds:
\begin{equation}\label{interp_Lagrange_formula}
Pf(x)=\sum\limits_{j=1}^{n+1}
f\left(x^{(j)}\right)\lambda_j(x). 
\end{equation}
Denote by $\|P\|_K$ the norm of $P$ as an operator from $C(K)$
in $C(K)$. From (\ref{interp_Lagrange_formula}), it~follows that
\begin{equation}\label{norm_of_P_lambda_j}
\|P\|_K=
\max_{x\in K}\sum_{j=1}^{n+1}
|\lambda_j(x)|.
\end{equation}
If $K$ is a convex polytope, then maximum in \eqref{norm_of_P_lambda_j} also can be taken  only over the~vertex set of  $K$. Although,  throughout the paper $K$ is a ball or an ellipsoid.

 By $\theta_n(K)$ we denote the minimal norm of an interpolation  projector
$P:C(K)\to \Pi_1({\mathbb R}^n)$ whose nodes belong to $K$.

Assume $K=B:=B(x^{(0)};R)$. It is proved in \cite{nevskii_2} that  
\begin{equation}\label{norm_P_ball_general}
\|P\|_B=
\max\limits_{f_j=\pm 1} \left[
R \left(\sum_{i=1}^n\left(\sum_{j=1}^{n+1} f_jl_{ij}\right)^2\right)^{1/2}
+\left|\sum_{j=1}^{n+1}f_j\lambda_j(x^{(0)})\right|
\right].
\end{equation}
If $S$ is a regular simplex inscribed into the ball, then $\|P\|_B$ does not depend neither on the center $x^{(0)}$, nor
on the radius $R$  of the ball, nor on the choice of such a simplex.
In this case,
$\|P\|_B=\max\{\psi(a_n),\psi(a_n+1)\}$ (see \cite[Theorem 2]{nevskii_2}). 
Here
\begin{equation}\label{psi_function_modulus}
\psi(t):=\frac{2\sqrt{n}}{n+1}\Bigl(t(n+1-t)\Bigr)^{1/2}+
\left|1-\frac{2t}{n+1}\right|, \quad  0\leq t\leq n+1,
\end{equation}
\begin{equation}\label{a_formula}
a_n:=\left\lfloor\frac{n+1}{2}-\frac{\sqrt{n+1}}{2}\right\rfloor.
\end{equation}


\section{The theorem on a simplex and its minimal  \\ ellipsoid}

Suppose $S$ is a nondegenerate $n$-dimensional simplex,   $E$  is a minimum volume ellipsoid containing $S$, and $m$ is
a natural number, $1\leq m\leq n$. To any set of $m$ vertices of $S$, let us put in correspondence a point $y\in E$
defined as follows. Let
$g$ \linebreak  be the center of gravity of the $(m-1)$-dimensional face of $S$ which contains the~chozen vertices, and let $h$ be the center of gravity of the $(n-m)$-dimensional face  which contains the rest  $n+1- m$ vertices. Then $y$ is the intersection point of the straight line $(gh)$ with the boundary of the ellipsoid in direction from $g$ to $h$.

\smallskip
{\bf Theorem 1.}  {\it For an arbitrary nondegenerate simplex $S\subset B_n$, there is such a~set ~of $m$  vertices for which $y\in B_n$.}

\smallskip 
{\it Proof.}  
Let $x^{(j)}$ be the vertices of $S$. The center of gravity of the simplex and also the center of gravity of its minimal
ellipsoid 
$E$ are  the point $c=\frac{1}{n+1}\sum x^{(j)}.$ Denote by $r$ the ratio of the distance between the center
of gravity of a regular simplex \linebreak and the~
center of gravity of its $(m-1)$-dimensional face to the radius of the circumscribed sphere. It is easy to see that
$$r=\frac{1}{m}\sqrt{m-\frac{m(m-1)}{n}}.$$   By $N$, we denote the number of $(m-1)$-dimensional
faces of an $n$-dimensional simplex: $N={n+1\choose m}$. 

Consider any set $J$ consisting of $m$ indices
$j.$ The center of gravity of the $(m-1)$-dimensional face with vertices  $x^{(j)}, j\in J,$ is the point
$g_J=\frac{1}{m}\sum_{j\in J} x^{(j)}$. Suppose $y_J$ coincides with the point $y$ for this set of vertices.
Then
$$y_J=c+\frac{1}{r}\left(c-g_J\right),$$
i.\,e., $y_J=(1/r)((r+1)c-g_J)$. 
Summing up over all sets $J$, we have
$$\sum_{J}\|y_J\|^2=\sum_{J}\langle y_J, y_J\rangle=$$
$$=\frac{1}{r^2}\left( N\|c\|^2(r+1)^2+\sum_J\|g_J\|^2-2(r+1)\langle
\sum_J g_J,c\rangle \right)=$$
\begin{equation}\label{right_side_of_sum_y_j_norms}
=\frac{1}{r^2}\left( N\|c\|^2(r^2-1)+\sum_J\|g_J\|^2 \right).
\end{equation}
We took into  account the equality $\sum g_J=Nc.$ Further, we claim that

\begin{equation}\label{norm_g_equality}
\sum_J\|g_J\|^2 = \frac{1}{mn}{n\choose m}\sum_{j=1}^{n+1} \|x^{(j)}\|^2+
\frac{(m-1)(n+1)}{mn}{n+1\choose m} \|c\|^2.
\end{equation}
To obtain \eqref{norm_g_equality}, let us remark  the following. The value $\sum \|g_J\|^2$ contains, for every
$j$, exactly ${n\choose m-1}$ numbers $\|x^{(j)}\|^2$ taken with coefficient $1/m^2$, and
exactly ${n-1\choose m-2}$  numbers
$2\langle x^{(i)},x^{(j)}\rangle$, for each $i\ne j$, with the same coefficient $1/m^2$. (The~latter ones exist for $m\geq 2$; in case $m=1$  there are no pairwise products.)  The expression $\|c\|^2$ contains all numbers $\|x^{(j)}\|^2$,  for any $j$,
and  numbers $2\langle x^{(i)},x^{(j)}\rangle$, for $i\ne j$, with~multiplier $1/(n+1)^2.$ The coefficient at
$\|c\|^2$ in the right-hand part of  \eqref{norm_g_equality} guarantees the~equality of~expressions with pairwise products $\langle x^{(i)},x^{(j)}\rangle$; the coefficient 
at~$\sum  \|x^{(j)}\|^2$ is chosen so that the difference of the left-hand part  of \eqref{norm_g_equality} 
and the second item  in~the right-hand part is equal to the first item. Since
$$ N(r^2-1)=-\frac{(m-1)(n+1)}{mn}{n+1\choose m},$$
after replacing $\sum \|g_J\|^2$ in 
\eqref{right_side_of_sum_y_j_norms} with the right-hand part of \eqref{norm_g_equality}, we notice that items with $\|c\|^2$ are compensated. After these  transformations,
we get 
\begin{equation}\label{means_equality}
\frac{1}{N}\sum_J\|y_J\|^2 =\frac{1}{n+1}\sum_{j=1}^{n+1} \|x^{(j)}\|^2.
\end{equation}
The inclusion $S\subset  B_n$ means that $\|x^{(j)}\|\leq 1.$ Therefore, the mean value of numbers $\|y_J\|^2$
is also $\leq 1.$ Consequently, for some set $J$ we have $\|y_J\|\leq 1,$
i.\,e., this point $y_J$ lies in $ B_n.$ Theorem 1 is proved.
\hfill$\Box$

\smallskip
As a conjecture, Theorem 1 was formulated by the author in \cite{nevskii_3}. 
Also in \cite{nevskii_3},  
the~statement of Theorem 1  in~case $m=1$ was proved 
in another way.

\section{The value $\theta_n(B_n)$} 

The main result of the paper concerns the minimum  norm of an ineterpolation projector with the nodes in an 
$n$-dimensional ball.

\smallskip
{\bf Theorem 2.} {\it  Let function $\psi(t)$ and number $a_n$ be defined by equalities 
\eqref{psi_function_modulus}--\eqref{a_formula}.
 For any $n\in{\mathbb N}$, 
 \begin{equation}\label{theta_equality}
 \theta_n(B_n)=\max\{\psi(a_n),\psi(a_n+1)\}.
 \end{equation}
 An interpolation projector 
  $P:C(B_n)\to \Pi_1({\mathbb R}^n)$  satisfies the equality $\|P\|_{B_n}=\theta_n(B_n)$ iff  the interpolation 
  nodes coincide with the vertices
  of a regular simplex inscribed into the boundary sphere.}

\smallskip 
{\it Proof.}  
First  let $P$ be an interpolation projector corresponding to a regular inscribed simplex  $S$ and let $\lambda_j(x)$ be
the basic Lagrange polynomials for this simplex. Denote
$p_n:=\|P\|_{B_n}$. It is proved in \cite{nevskii_2} that $p_n=\max\{\psi(a_n),\psi(a_n+1)\}.$ 
The points $y\in B_n$ satisfying
\begin{equation}\label{p_n_equality}
\sum_{j=1}^{n+1} |\lambda_j(y)|=p_n
\end{equation}
 were found in \cite{nevskii_3}.
Namely, assume $k_n$ coincides with that of the numbers $a_n$ and $a_n+1$ which delivers to $\psi(t)$ the maximum value. Let us fix $m=k_n$. Then equality \eqref{p_n_equality} takes place for every point   $y$
which is constructed above Theorem 1 in case  $E=B_n$. The number of these points 
$y$ is $N={n+1\choose m}$.

Now suppose $P:C(B_n)\to \Pi_1({\mathbb R}^n)$ is an arbitrary projector with the nodes
$x^{(j)}\in B_n$, $S$ is the simplex with  the vertices $x^{(j)}$, and $E$ is a minimum volume ellipsoid
contained $S$. We can also consider $P$ as a projector acting from $C(E)$.
Then
$$\|P\|_E=p_n=\sum_{j=1}^{n+1} |\lambda_j(y)|$$
for $m=k_n$ and for each of the points $y$ noted at the beginning of the previous section (now with respect to the minimal ellipsoid
$E$).
By Theorem 1, at least one of these points, say $y^*$, belongs to $B_n$. Therefore,
$$\|P\|_{B_n}=\max_{x\in B_n} \sum_{j=1}^{n+1}|\lambda_j(x)|\geq\sum_{j=1}^{n+1}
 |\lambda_j(y^*)|=p_n.$$
 Thus, $\theta_n(B_n)=p_n=\max\{\psi(a_n),\psi(a_n+1)\}=\psi(k_n)$.

Let us demonstrate that if $\|P\|_{B_n}=\theta_n(B_n)$, then  the simplex with vertices at the interpolation nodes
is  inscribed into $B_n$ and regular. For this simplex $S$, some point  $y\in E$,
constructed for $m=k_n$, falls on the boundary sphere --- otherwise we have
$\|P\|_{B_n}>p_n=\theta_n(B_n)$. This is due to the fact that $\sum |\lambda_j(x)|$ 
increases  monotonously as $x$ moves
in a straight line in direction from $c$ to $y$. Since $\|y\|=1,$ the mean value of $\|y_J\|$ for the rest
$y_J$ is also $\leq 1$ (see \eqref{means_equality}). Hence, another such point also lies in the boundary of the ball. Thus we obtain $\|y_J\|=1$ for all sets $J$ consisting of $m=k_n$ indices $j$. Now \eqref{means_equality} yields 
$\|x^{(j)}\|=1,$ i.\,e., the simplex is inscribed into the ball. Since all the points
$y_J$ belong to the boundary sphere, the function $\sum |\lambda_j(x)|$  has maximum upon the ball at
$N={n+1\choose m}$ different points. Consequently, the simplex  $S$ is regular.
This completes the proof.
\hfill$\Box$

\smallskip
For $1\leq n\leq 4$, equality \eqref{theta_equality} and characterization of minimal projectors were obtained in \cite{nevskii_2} and \cite{nevskii_3} by another methods.

\smallskip 
{\bf Corollary.} {\it Always
$\sqrt{n}\leq \theta_n(B_n)\leq \sqrt{n+1}.$
Moreover, $\theta_n(B_n)$ $=$ $\sqrt{n}$ only for~$n=1$ and
$\theta_n(B_n)=\sqrt{n+1}$ when and only when
$\sqrt{n+1}$  is an integer.}

\smallskip
 Really, it was proved in  \cite{nevskii_2} that the above relations are satisfied for  values
 $p_n=\max\{\psi(a_n),\psi(a_n+1)\}$. It remains to take into account Theorem 2.

\smallskip
Note that the numbers $k_n$ from the proof of Theorem 2 increase with  $n$, but not strictly  monotonously.
If $n\geq 2$, then
$k_n\leq n/2.$ 
The numbers $a_n, k_n $, and~$ \theta_n(B_n)=\psi(k_n)$ are given in Table 1; we use the results from~\cite{nevskii_2}.

\begin{table}[h!!]
\begin{center}
\caption{The numbers $a_n, k_n $, and $ \theta_n(B_n)=\psi(k_n)
$}
\label{tab:nev_uhl_norm_est}
\bigskip
\bgroup
$
\def\arraystretch{1.7}
\begin{array}{|c|c|c|c|c|c|c|}
\hline
 n & a_n & a_n+1 & \psi(a_n) & \psi(a_n+1) & k_n&  \theta_n(B_n)=\psi(k_n)
\\
\hline
 1 & 0 & 1 & 1 & 1 & 1 & 1 \\
 \hline
 2 & 0 & 1 & 1 & \frac{5}{3} & 1 & \frac{5}{3} =1.6666\ldots \\
 \hline
 3 & 1 & 2 & 2 & \sqrt{3} & 1 & 2 \\
 \hline
 4 & 1 & 2 & \frac{11}{5} & \frac{1+4 \sqrt{6}}{5}  & 1 & \frac{11}{5} = 2.2 \\
 \hline
 5 & 1 & 2 & \frac{7}{3} & \frac{1+2 \sqrt{10}}{3}  & 2 & \frac{1+2 \sqrt{10}}{3} =2.4415\ldots \\
 \hline
 6 & 2 & 3 & \frac{3+4 \sqrt{15}}{7}  & \frac{ 1+12 \sqrt{2}}{7} & 2 & 
 \frac{3+4 \sqrt{15}}{7}=2.6417\ldots \\
 \hline
 7 & 2 & 3 & \frac{1+\sqrt{21}}{2}  & \frac{1+\sqrt{105}}{4}  & 3 & \frac{1+\sqrt{105}}{4}  = 2.8117\ldots \\
 \hline
 8 & 3 & 4 & 3 & \frac{1+8 \sqrt{10}}{9}  & 3 & 3 \\
 \hline
 9 & 3 & 4 & \frac{2+3 \sqrt{21}}{5}  & \frac{1+6 \sqrt{6}}{5}  & 3 & \frac{2+3 \sqrt{21}}{5}  = 3.1495\ldots \\
 \hline
 10 & 3 & 4 & \frac{ 5+8 \sqrt{15}}{11} & \frac{ 3+4 \sqrt{70}}{11} & 4 & 
 \frac{3+4 \sqrt{70}}{11}  = 3.3151\ldots \\
 \hline
 11 & 4 & 5 & \frac{ 1+2 \sqrt{22}}{3} & \frac{ 1+\sqrt{385}}{6} & 4 & \frac{1+2 \sqrt{22}}{3}  = 3.4602\ldots \\
 \hline
 12 & 4 & 5 & \frac{5+24 \sqrt{3}}{13}  & \frac{3+8 \sqrt{30}}{13}  & 5 & 
 \frac{3+8 \sqrt{30}}{13}  = 3.6013\ldots \\
 \hline
 13 & 5 & 6 & \frac{2+3 \sqrt{65}}{7}  & \frac{1+4 \sqrt{39}}{7}  & 5 & 
 \frac{2+3 \sqrt{65}}{7}  = 3.7409\ldots \\
 \hline
 14 & 5 & 6 & \frac{1+4 \sqrt{7}}{3}  & \frac{1+4 \sqrt{21}}{5}  & 6 & 
 \frac{1+4 \sqrt{21}}{5}  = 3.8660\ldots \\
 \hline
 15 & 6 & 7 & 4 & \frac{1+3 \sqrt{105}}{8}  & 6 & 4  \\
 \hline
 50 & 21 & 22 & \frac{3+20 \sqrt{35}}{17}  & \frac{7+20 \sqrt{319}}{51}  & 22 & 
 \frac{7+20 \sqrt{319}}{51}  = 7.1414\ldots \\
 \hline
 100 & 45 & 46 & \frac{11+120 \sqrt{70}}{101}  & \frac{ 9+20 \sqrt{2530}}{101} & 45 & \frac{11+120 \sqrt{70}}{101}  = 10.0494\ldots \\
 \hline
 1000  & 484 & 485 & \frac{3+40 \sqrt{5170}}{91}  & \frac{31+200 \sqrt{25026}}{1001} & 485 & \frac{31+200 \sqrt{25026}}{1001} = 31.6385\ldots \\
 \hline
\end{array}
$
\egroup
\end{center}
\end{table}

\smallskip
\section{Acknowledgement}
The author is grateful to Arseniy Akopyan for his help while proving Theorem~1.


\end{document}